# G ADD-ON, DIGITAL, SIEVE, GENERAL PERIODICAL, AND NON-ARITHMETIC SEQUENCES

by Florentin Smarandache, Ph. D.
University of New Mexico
Gallup, NM 87301, USA

**Abstract**:
Other new sequences are introduced in number theory, and for each one a general question: how many primes each sequence has.

**Keywords**: sequence, symmetry, consecutive, prime, representation of numbers.

**1991 MSC**: 11A67

**Introduction.**
In this paper is presented a small survey on fourteen sequences, such as: G Add-on Sequences, Sieve Sequences, Digital Sequences, Non-Arithmetic Sequences, recreational sequences (Lucky Method/Operation/Algorithm/Differentiation/Integration etc.), General Periodical Sequences, and arithmetic functions.

**1) G Add-On Sequence (I)**

Let $G = \{g_1, g_2, \ldots, g_k, \ldots\}$ be an ordered set of positive integers with a given property G.
Then the corresponding G Add-On Sequence is defined through

$$SG = \{a_i : a_1 = g_1, a_k = a_{k-1} \cdot 10^{1+\log_{10}(g_k)} + g_k, k \geq 1\}.$$

H. Ibstedt studied some particular cases of this sequence, that he has presented to the FIRST INTERNATIONAL CONFERENCE ON SMARANDACHE TYPE NOTIONS IN NUMBER THEORY, University of Craiova, Romania, August 21-24, 1997.

*a) Examples of G Add-On Sequences (II)*

The following particular cases were studied:
  *a.1) Odd Sequence is generated by choosing*

G = {1, 3, 5, 7, 9, 11, ...}, and it is:
1, 13, 135, 1357, 13579, 1357911, 13571113, ... .
Using the elliptic curve prime factorization program we find
the first five prime numbers among the first 200 terms of this
sequence, i.e. the ranks 2, 15, 27, 63, 93.
But are they infinitely or finitely many?

*a.2) Even Sequence is generated by choosing*
G = {2, 4, 6, 8, 10, 12, ...}, and it is:
2, 24, 246, 2468, 246810, 24681012, ... .
Searching the first 200 terms of the sequence we didn't find
any n-th perfect power among them, no perfect square, nor even
of the form 2p, where p is a prime or pseudo-prime.
Conjecture: There is no n-th perfect power term!

*a.3) Prime Sequence is generated by choosing*
G = {2, 3, 5, 7, 11, 13, 17, ...}, and it is:
2, 23, 235, 2357, 235711, 23571113, 2357111317, ... .
Terms #2 and #4 are primes; terms #128 (of 355 digits) and #174
(of 499 digits) might be, but we couldn't check -- among the first
200 terms of the sequence.
Question: Are there infinitely or finitely many such primes?
(H. Ibstedt)

Reference:
[1] Mudge, Mike, "Smarandache Sequences and Related Open Problems",
    Personal Computer World, Numbers Count, February 1997.

## 2) Non-Arithmetic Progressions (I)

One of them defines the t-Term Non-Arithmetic Progression
as the set:
$\{a_i : a_i \text{ is the smallest integer such that } a_i > a_{i-1},$
and there are at most t-1 terms in an arithmetic progression$\}$.
A QBASIC program was designed to implement a strategy for
building a such progression, and a table for the 65 first terms of
the non-arithmetic progressions for t=3 to 15 is given
(H. Ibstedt).

Reference:
[1] Mudge, Mike, "Smarandache Sequences and Related Open Problems",
    Personal Computer World, Numbers Count, February 1997.

## 3) Concatenation Type Sequences

Let $s_1, s_2, s_3, ..., s_n, ...$ be an infinite integer sequence
(noted by S). Then the Concatenation is defined as:

$s_1, \overline{s_1 s_2}, \overline{s_1 s_2 s_3}, ...$ .

H. Ibstedt searched, in some particular cases, how many terms of this
concatenated S-sequence belong to the initial S-sequence.

**4) Partition Type Sequences**

   Let f be an arithmetic function, and R a k-relation among numbers.
   How many times can n be expressed under the form of;
   $n = R(f(n_1), f(n_2), ..., f(n_k))$,

   for some k and $n_1, n_2, ..., n_k$ such that $n_1 + n_2 + ... + n_k = n$ ?

 Look at some particular cases:  How many times can be n express as a
 sum of non-null squares (or cubes, or m-powers)?

**5) The Lucky Method/Algorithm/Operation/Differentiation/Integration/etc.**

   Generally:  The Lucky
Method/Algorithm/Operation/Differentiation/Integration/etc. is
said to be any incorrect method or algorithm or operation etc. which leads to
a correct result.  The wrong calculation should be fun, somehow similarly
to the students' common mistakes, or to produce confusions or paradoxes.
   Can someone give an example of a Lucky Derivation, or
Integration, or Lucky Solution to a Differential Equation?

As a particular case:
A number is said to be an *S. Lucky Number* if an incorrect
calculation leads to a correct result, which is that number.
Is the set of all fractions, where such (or another) incorrect calculation
leads to a correct result, finite or infinite?

Reference:
  [1] Smarandache, Florentin, "Collected Papers" (Vol. II), University of
      Kishinev, 1997.

**6)**   Construction of Elements of the **Square-Partial-Digital Subsequence**

The Square-Partial-Digital Subsequence (SPDS) is the
sequence of square integers which admit a partition for which each segment is
a square integer.  An example is $506^2 = 256036$, which has partition 256/0/36.
C. Ashbacher showed that SPDS is infinite by exhibiting two infinite families
of elements.  We will extend his results by showing how to construct infinite
families of elements of SSPDS containing desired patterns of digits.
   Unsolved Question 1:
441 belongs to SSPDS, and his square $441^2 = 194481$ also belongs to SSPDS.
Can an example be found of integers $m, m^2, m^4$ all belonging to SSPDS?
   Unsolved Question 2:

It is relatively easy to find two consecutive squares in SSDPS, i.e.
$12^2 = 144$ and $13^2 = 169$.
Does SSDPS also contain three or more consecutive squares?
What is the maximum length?

**7)  Prime-Digital Sub-Sequence**

"Personal Computer World" Numbers Count of February 1997
presented some of the Smarandache Sequences and related open
problems.
One of them defines the Prime-Digital Sub-Sequence
as the ordered set of primes whose digits are all primes:
2, 3, 5, 7, 23, 37, 53, 73, 223, 227, 233, 257, 277, ... .
H. Ibstedt used a computer program in Ubasic to calculate the first 100
terms of the sequence.  The 100-th term is 33223.
Sylvester Smith [2] conjectured that the sequence is infinite.  We
also agree that this sequence is in fact infinite.

References:
   [1] Mudge, Mike, "Smarandache Sequences and Related Open Problems",
       Personal Computer World, Numbers Count, February 1997.
   [2] Smith, Sylvester, "A Set of Conjectures on Smarandache
       Sequences", in <Bulletin of Pure and Applied Sciences>,
       Delhi, India, Vol. 15E, No. 1, 1996, pp. 101-107.

**8) Special Expressions**

*a) Perfect Powers in Special Expressions (I)*

How many primes are there in the Special Expression:
$$x^y + y^x,$$
where $\gcd(x, y) = 1$ ? [J. Castillo & P. Castini]
K. Kashihara announced that there are only finitely many numbers of the
above form which are products of factorials.
F. Luca proposed the following conjecture:
  Let a, b, and c three integers with ab nonzero.  Then the equation:
  $ax^y + by^x = cz^n$, with $x, y, n \geq 2$, and $\gcd(x, y) = 1$,
  has finitely many solutions (x, y, z, n).

*b)  Products of Factorials in Smarandache Type Expressions (II)*

J. Castillo ["Mathematical Spectrum", Vol. 29, 1997/8, 21] asked
how many primes are there in the Smarandache n-Expression:
$$x_1^{x_2} + x_2^{x_3} + \ldots + x_n^{x_1},$$
where $n > 1$, $x_1, x_2, \ldots, x_n > 1$, and $\gcd(x_1, x_2, \ldots, x_n) = 1$ ?
[This is a generalization of the Smarandache 2-Expression:  $x^y + y^x$.]

F. Luca announced a lower bound for the size of the largest prime
   divisor of an expression of type $ax^y + by^x$, where $ab$ is nonzero,
   $x, y \geq 2$, and $\gcd(x, y) = 1$.

## 9) The General Periodic Sequence

Definition:
Let S be a finite set, and $f : S \longrightarrow S$ be a function defined
for all elements of S.
There will always be a periodic sequence whenever we repeat the composition
of the function f with itself more times than card(S), accordingly to the
box principle of Dirichlet.
[The invariant sequence is considered a periodic sequence whose period length
has one term.]

Thus the General Periodic Sequence is defined as:
  $a_1 = f(s)$, where s is an element of S;
  $a_2 = f(a_1) = f(f(s))$;
  $a_3 = f(a_2) = f(f(a_1)) = f(f(f(s)))$;
  and so on.

M. R. Popov particularized S and f to study interesting cases of this type of
sequences.

## 10) n-Digit Periodical Sequences

*a) The Two-Digit Periodic Sequence (I)*

 Let $N_1$ be an integer of at most two digits and let $N_1'$ be its
 digital reverse. One defines the absolute value $N_2 = \mathrm{abs}(N_1 - N_1')$.
 And so on: $N_3 = \mathrm{abs}(N_2 - N_2')$, etc. If a number N has one digit only,
 one considers its reverse as $N \times 10$ (for example: 5, which is 05, reversed will
 be 50). This sequence is periodic.
 Except the case when the two digits are equal, and the sequence becomes:
   $N_1, 0, 0, 0, \ldots$
 the iteration always produces a loop of length 5, which starts on the second
 or the third term of the sequence, and the period is 9, 81, 63, 27, 45
 or a cyclic permutation thereof.

 Reference:
   [1] Popov, M.R., "Smarandache's Periodic Sequences", in <Mathematical
   Spectrum>, University of Sheffield, U.K., Vol. 29, No. 1, 1996/7, p. 15.

(The next periodic sequences are extracted from this paper too).

*b) The n-Digit Smarandache Periodic Sequence (II)*

Let $N_1$ be an integer of at most n digits and let $N_1'$ be its
 digital reverse. One defines the absolute value $N_2 = \mathrm{abs}(N_1 - N_1')$.

And so on: N3 = abs (N2 - N2'), etc.  If a number N has less than n digits,
   one considers its reverse as N'x(10^k), where N' is the reverse of N and
   k is the number of missing digits, (for example: the number 24 doesn't have
   five digits, but can be written as 00024, and reversed will be 42000).
   This sequence is periodic according to Dirichlet's box principle.
      The Smarandache 3-Digit Periodic Sequence (domain 100 <= N1 <= 999):
   - there are 90 symmetric integers, 101, 111, 121, ..., for which N2 = 0;
   - all other initial integers iterate into various entry points of the same
   periodic subsequence (or a cyclic permutation thereof) of five terms:
                     99, 891, 693, 297, 495.
      The Smarandache 4-Digit Periodic Sequence (domain 1000<= N1 <= 9999):
   - the largest number of iterations carried out in order to reach the first
   member of the loop is 18, and it happens for N1 = 1019;
   - iterations of 8818 integers result in one of the following loops (or a
   cyclic permutation thereof): 2178, 6534;  or 90, 810, 630, 270, 450;  or
   909, 8181, 6363, 2727, 4545;  or 999, 8991, 6993, 2997, 4995;
   - the other iterations ended up in the invariant 0.
   (H. Ibstedt)

*c)  The 5-Digit and 6-Digit Smarandache Periodic Sequences (III)*

Let N1 be an integer of at most n digits and let N1' be its
 digital reverse.  One defines the absolute value N2 = abs (N1 - N1').
 And so on: N3 = abs (N2 - N2'), etc.  If a number N has less than n digits,
 one considers its reverse as N'x(10^k), where N' is the reverse of N and
 k is the number of missing digits, (for example: the number 24 doesn't have
 five digits, but can be written as 00024, and reversed will be 42000).
 This sequence is periodic according to Dirichlet's box principle, leading to
 invariant or a loop.
    The Smarandache 5-Digit Periodic Sequence (domain 10000 <= N1 <= 99999):
 - there are 920 integers iterating into the invariant 0 due to symmetries;
 - the other ones iterate into one of the following loops (or a cyclic
 permutation of these): 21978, 65934;  or 990, 8910, 6930, 2970, 4950;  or
 9009, 81081, 63063, 27027, 45045;  or 9999, 89991, 69993, 29997, 49995.
    The Smarandache 6-Digit Periodic Sequence (domain 100000 <= N1 <= 999999):
 - there are 13667 integers iterating into the invariant 0 due to symmetries;
 - the longest sequence of iterations before arriving at the first loop
 member is 53 for N1 = 100720;
 - the loops have 2, 5, 9, or 18 terms.

*d)  The Subtraction Periodic Sequences (IV)*

Let c be a positive integer.  Start with a positive integer N, and
let N' be its digital reverse.  Put N1 = abs(N1' - c), and let N1' be its
digital reverse.  Put N2 = abs (N1' - c), and let N2' be its digital reverse.
And so on.  We shall eventually obtain a repetition.
For example, with c = 1 and N = 52 we obtain the sequence:  52, 24, 41, 13,
30, 02, 19, 90, 08, 79, 96, 68, 85, 57, 74, 46, 63, 35, 52, ... .  Here a
repetition occurs after 18 steps, and the length of the repeating cycle is 18.
   First example:  c = 1, 10<= N <= 999.
Every other member of this interval is an entry point into one of five cyclic
periodic sequences (four of these are of length 18, and one of length 9).
When N is of the form 11k or 11k-1, then the iteration process results in 0.

Second example:  1 <= c <= 9, 100 <= N <= 999.
For c = 1, 2, or 5 all iterations result in the invariant 0 after, sometimes,
a large number of iterations.
For the other values of c there are only eight different possible values for
the length of the loops, namely 11, 22, 33, 50, 100, 167, 189, 200.
For c = 7 and N = 109 we have an example of the longest loop obtained:  it
has 200 elements, and the loop is closed after 286 iterations.
(H. Ibstedt)

e)  *The Multiplication Periodic Sequences (V)*

Let c > 1 be a positive integer.  Start with a positive integer N,
multiply each digit x of N by c and replace that digit by the last digit of
cx to give N1.  And so on.  We shall eventually obtain a repetition.
For example, with c = 7 and N = 68 we obtain the sequence:
     68, 26, 42, 84, 68, ... .
Integers whose digits are all equal to 5 are invariant under the given
operation after one iteration.

One studies the One-Digit Multiplication Periodic Sequences only.
(For c of two or more digits the problem becomes more complicated.)
  If c = 2, there are four term loops, starting on the first or second term.
  If c = 3, there are four term loops, starting with the first term.
  If c = 4, there are two term loops, starting on the first or second term
(could be called Switch or Pendulum Sequence).
  If c = 5 or 6, the sequence is invariant after one iteration.
  If c = 7, there are four term loops, starting with the first term.
  If c = 8, there are four term loops, starting with the second term.
  If c = 9, there are two term loops, starting with the first term (pendulum).
(H. Ibstedt)

e)  *The Mixed Composition Periodic Sequences (VI)*

Let N be a two-digit number.  Add the digits, and add them again if
the sum is greater than 10.  Also take the absolute value of their difference.
These are the first and second digits of N1.  Now repeat this.
For example, with N = 75 we obtain the sequence:  75, 32, 51, 64, 12, 31, 42,
62, 84, 34, 71, 86, 52, 73, 14, 53, 82, 16, 75, ... .
There are no invariants in this case.  Four numbers: 36, 90, 93, and 99
produce two-element loops.  The longest loops have 18 elements.  There also
are loops of 4, 6, and 12 elements.
(H. Ibstedt)

   There will always be a periodic (invariant) sequence whenever we have a
function f : S ---> S, where S is a finite set,
and we repeat the function f more times than card(S).
Thus the General Periodic Sequence is defined as:
  a1 = f(s), where s is an element of S;
  a2 = f(a1) = f(f(s));
  a3 = f(a2) = f(f(a1)) = f(f(f(s)));
  and so on.

**11) New Sequences: The Family of Metallic Means**

The family of Metallic Means (whom most prominent
members are the Golden Mean, Silver Mean, Bronze Mean, Nickel Mean, Copper
Mean, etc.) comprises every quadratic irrational number that is the
positive solution of one of the algebraic equations
$$x^2 - nx - 1 = 0 \quad \text{or} \quad x^2 - x - n = 0,$$
where n is a natural number.
All of them are closely related to quasi-periodic dynamics, being therefore
important basis of musical and architectural proportions.  Through the
analysis of their common mathematical properties, it becomes evident that
they interconnect different human fields of knowledge, in the sense defined
by Florentin Smarandache ("Paradoxist Mathematics").
Being irrational numbers, in applications to different scientific
disciplines, they have to be approximated by ratios of integers -- which is
the goal of this paper.
(Vera W. de Spinadel)

**12)  Some New Functions in the Number Theory**

We investigate and prove the functions:
$S_1 : N-\{0,1\} \longrightarrow N$, $S_1(n) = 1/S(n)$;
$S_2 : N^* \longrightarrow N$, $S_2(n) = S(n)/n$
verify the Lipschitz condition, but the functions:
$S_3 : N-\{0,1\} \longrightarrow N$, $S_3(n) = n/S(n)$;
$F_s : N^* \longrightarrow N$,
$$F_s(x) = \sigma( S(p_i^x) ) \text{ for } i \text{ from } 1 \text{ to } \pi(x),$$
where $p_i$ are the prime numbers not greater than x and
$\pi(x)$ is the number of them;
$\Theta : N^* \longrightarrow N$,
$$\Theta(x) = \sigma\, S(p_i^x), \text{ where } p_i \text{ are prime numbers}$$
which divide x;
$\overline{\Theta} : N^* \longrightarrow N$,
$$\overline{\Theta}(x) = \sigma\, S(p_i^x), \text{ where } p_i \text{ are prime numbers}$$
which do not divide x;
where $S(n)$ is the Smarandache function for all six previous functions,
verify the Lipschitz condition.
(V. Seleacu, S. Zanfir)

Reference:
  [1] Mencze, M., "Smarandache Relationships and Subsequences", <Bulletin
      of Pure and Applied Sciences>, Delhi, India, Vol. 17E, No. 1, pp.
      55-62, 1998.

**13) Erdos-Smarandache Numbers:**

2, 3, 5, 6, 7, 10, 11, 13, 14, 15, 17, 19, 20, 21, 22, 23, 26, 28, 29,
    30, 31, 33, 34, 35, ... .
    Solutions to the diophantine equation P(n)=S(n), where P(n) is the
    largest prime factor which divides n, and S(n) is the classical
    Smarandache function: the smallest integer such that S(n)! is
    a multiple of n.

References:
   [1] Erdos, P., Ashbacher C., Thoughts of Pal Erdos on Some Smarandache
       Notions, <Smarandache Notions Journal>, Vol. 8, No. 1-2-3, 1997,
       220-224.
   [2] Sloane, N. J. A., On-Line Encyclopedia of Integers, Sequence
       A048839.

**14) n-ary Sieve**:
    1, 2, 4, 7, 9, 14, 20, 25, 31, 34, 44, ... .
    Keep the first k numbers, skip the k+1 numbers, for k = 2, 3, 4, ... .

References:
   [1] Le M., "On the Smarandache n-ary sieve", <Smarandache Notions Journal>,
       Vol. 10, No. 1-2-3, 1999, 146-147.
   [2] Sloane, N. J. A., On-Line Encyclopedia of Integers, Sequences
       A048859, A007952.

           **General References** (edited by C. Dumitrescu & M. Bencze):

   1.   Sloane, N. J. A., Plouffe, Simon, "The On-Line Encyclopedia of
        Integer Sequences", AT&T Bell Labs, Murray Hill, New Jersey, USA,
        1995, online: superseeker@research.att.com ;
   2.   Mudge, Mike, "Some Sequences of Smarandache contrasted with a
        permutations problem relating to class/student allocation", in
        <Personal Computer World>, London, England, June 1995, pp. 674-5;
   3.   Lohon, O., Buz, Maria, University of Craiova Library, Letter No.
        499, July 07, 1995.
   4.   Mudge, Mike, "Pseudo-Sequences (according to Florentin  Smarandache)
        and Generalized Russian Multiplication", in <Personal Computer
        World>, London, England, August 1995, pp. 671;
   5.   Growney, JoAnne, Bloomsburg University, PA, "The most Humanistic
        Mathematician: Florentin Smarandache"
        and Larry Seagull, "Poem in Arithmetic Space", in the <Humanistic
        Mathematics Network>, Harvey Mudd College, Claremont, CA, October
        1995, # 12, p. 38 and pp. 38-40 respectively;
   6.   Le, Charles T., "The most paradoxist mathematician of the world", in
        <Bulletin of Number Theory>, March 1995, Vol. 3, No. 1;
   7.   Seagull, Larry, Glendale Community College, "Poem in Arithmetic
        Space", in <Abracadabra>, Salinas, CA, August 1995, Anul III, No.
        34, pp.20-1.
   8.   Moore, Carol (Library Specialist), Wurzburger, Marilyn (Head of
        Specila Collections), Abstract of "The Florentin Smarandache papers"
        special collection, Call # SM SC SM-15, at Arizona State University,
        Tempe, AZ 85287-1006, Box 871006, Tel. (602) 965-6515, E-mail:
        icclmc@asuvm.inre.asu.edu , USA;